\documentclass[a4paper, 10pt, twoside]{amsart}

\usepackage[centertags]{amsmath}
\usepackage{amsthm}
\usepackage[DIV9, headinclude, BCOR 10mm]{typearea}
\usepackage{omega_math}
\usepackage{color}
\usepackage{tabularx}
\usepackage{url}
\usepackage{omega_layout}
\usepackage{calc}
\usepackage{ifpdf}
\ifpdf
  \usepackage[all]{xy}
  \SelectTips{cm}{}
  \usepackage[pdftex]{graphicx}
  \usepackage[pdftex,pagebackref,colorlinks,linkcolor=blue,citecolor=blue,urlcolor=blue,a4paper,hypertexnames=false]{hyperref}
  \hypersetup{pdfauthor=Clara L\"oh,pdftitle=Fundamental classes of
    negatively curved manifolds cannot be represented by products of surfaces}
  
\else
  \usepackage[all, dvips]{xy}
  \SelectTips{cm}{}
  \usepackage{graphicx}
  
\fi
\usepackage{omega_thm} 
\usepackage{caption}

\title[Fundamental classes not representable by products]{Fundamental
  classes\\ of negatively curved manifolds\\ cannot be
  represented by products of manifolds}
\author{Clara L\"oh}
\date{\today}

\begin{document}

\begin{abstract}
  Not every singular homology class is the push-forward of the
  fundamental class of some manifold. In the same spirit, one can
  study the following problem: Which singular homology classes are the
  push-forward of the fundamental class of a given type of manifolds?
  In the present article, we show that the fundamental classes of
  negatively curved manifolds cannot be represented by a non-trivial
  product of manifolds. This observation sheds some light on the
  functorial semi-norm on singular homology given by products of
  compact surfaces.
\end{abstract}

\maketitle
\thispagestyle{empty}

\section{Introduction}

\noindent
As manifolds form a particularly accessible class of topological
spaces, it is natural to ask which classes in singular homology with
$\Z$-coefficients are the push-forward of the fundamental class of
some manifold; indeed, not every singular homology class is such a
push-forward~\cite[Th\'eor\`eme~3]{thom}. In the same spirit, one can
study the following more restrictive problem: Which singular homology
classes are the push-forward of the fundamental class of a given type
of manifolds? An interesting type of manifolds in the context of
functorial semi-norms on singular homology~\cite[5.34
on~p.~302]{gromovmetric} is products of compact surfaces.

\begin{defi}
  Let $M$ be an oriented, closed, connected manifold. We say that the
  fundamental class of~$M$ can be \emph{represented by a product of
    manifolds}, if there exist~$d \in \Z\setminus\{0\}$,~~$r
  \in \N_{>1}$, and oriented, closed, connected manifolds~$S_1, \dots, S_r$
  of non-zero dimension admitting a continuous map~$f \colon S_1
  \times \dots \times S_r \longrightarrow M$ such that
  \[ \hsing * {f; \Z}\bigl( \fcl{S_1 \times \dots \times S_r}\bigr)
     = d \cdot \fcl M.
     \qedhere
  \]
\end{defi}

Gromov suspected that most interesting homology classes -- such as
fundamental classes of irreducible locally symmetric spaces -- cannot
be represented by products of manifolds~\cite[5.36
on~p.~303f]{gromovmetric}; in the present article, we confirm this
anticipation, at least in the case of manifolds of negative curvature:

\begin{satz}\label{hypthm}
  Let $M$ be an oriented, closed, connected, smooth manifold that
  admits a Riemannian metric with negative sectional curvature. Then
  the fundamental class of~$M$ cannot be represented by a product of
  manifolds, let alone by a product of surfaces.
\end{satz}

In particular, we obtain that the functorial semi-norm given by
products of compact surfaces~\cite[5.36 on~p.~303f]{gromovmetric} is
infinite when evaluated on the fundamental class of an oriented,
closed, connected Riemannian manifold with negative sectional
curvature.

In order to prove Theorem~\ref{hypthm}, we assume that the fundamental
class of~$M$ can be represented by a product of manifolds and lead
this assumption to a contradiction by studying the corresponding scene
on the level of fundamental groups. Roughly speaking, the source of
this contradiction lies in the tension between $\pi_1(M)$ being a
complicated group (because $M$ is negatively curved) and $\pi_1(M)$
being too commutative (the images of the fundamental groups of the
factors commute in~$\pi_1(M)$). 

In fact, we present two directions of generalisations of
Theorem~\ref{hypthm} in Sections~\ref{asphericalsubsec}
and~\ref{simvolsubsec} below, and establish Theorem~\ref{hypthm} as an
instance of these more general statements.

\subsection{Irrepresentability -- aspherical manifolds}\label{asphericalsubsec}

The first generalisation applies to a larger class of aspherical
manifolds and provides an obstruction to representability by a product
of manifolds in terms of group theory; in order to formulate the
statement, we introduce the following conventions:

\begin{defi}\label{classdef}\hfil
  \begin{itemize}
    \item Let $C$ be a class of groups. We say that a group~$G$
      \emph{has centralisers in}~$C$, if for every element~$g\in G$ of
      infinite order the centraliser~$C_G(g)$ lies in~$C$.  
    \item A class of groups is a \emph{class with products}, if it is
      closed under taking subgroups, quotients, and finite products.
    \qedhere
  \end{itemize}
\end{defi}

\begin{satz}\label{asphthm}
  Let $M$ be an aspherical, oriented, closed, connected manifold whose
  fundamental class can be represented by a product of manifolds. If
  the fundamental group of~$M$ has centralisers in a class~$C$ of
  groups with products, then the fundamental group of~$M$ contains a
  subgroup of finite index lying in~$C$ (i.e., $\pi_1(M)$ is virtually
  in~$C$).
\end{satz}

\subsection{Irrepresentability -- manifolds with positive simplicial volume}\label{simvolsubsec}

For the second generalisation, we trade the asphericity condition
for a suitable hypothesis on the simplicial volume, an invariant
unaware of higher homotopy.

\begin{satz}\label{svthm}
  Let $M$ be an oriented, closed, connected manifold with non-zero
  simplicial volume whose fundamental group is large and has amenable
  centralisers. Then the fundamental class of~$M$ cannot be
  represented by a product of manifolds.
\end{satz}

We now explain the occurring terminology in more detail:

\begin{defi}\label{largedef}\hfil
   A group is \emph{large}, if it contains an element of infinite order.
\end{defi}

Recall that a group~$G$ is called \emph{amenable} if there is a
left-invariant mean on the set~$B(G,\R)$ of bounded functions from~$G$
to~$\R$, i.e., if there is a $G$-invariant, linear map~$m \colon B(G,\R)
\longrightarrow \R$ with~$\inf_{g \in G} f(g) \leq m(f) \leq
\sup_{g\in G} f(g)$ for all~$f \in B(G,\R)$. 

Every finite, every Abelian, and every solvable group is amenable. The
class of amenable groups is closed with respect to taking subgroups,
quotients, and extensions; in particular, the class of amenable groups
is a class with products in the sense of Definition~\ref{classdef}.
The group~$\Z * \Z$ is not amenable. A thorough account of amenable
groups is given in Paterson's book~\cite{paterson}.

For a topological space~$X$ and a class~$\alpha \in \hsing k X$ in
singular homology with \mbox{$\R$-co}\-ef\-ficients, the $\ell^1$\emph{-semi-norm
  of}~$\alpha$ is defined by
  \[ \lone \alpha := \inf \bigl\{ \lone c
                          \bigm| \text{$c \in \csing k X$ is an $\R$-cycle representing~$\alpha$}
                          \bigr\};
  \]
  here, for a singular chain~$c = \sum_{j=0}^r a_j \cdot \sigma_j \in
  \csing k X$, we write~$\lone c := \sum_{j=0}^r |a_j|$.

  The simplicial volume is an example of a topological invariant
  defined in terms of the $\ell^1$-semi-norm~\cite[p.~8]{vbc}: The
  \emph{simplicial volume} of an oriented, closed, connected
  manifold~$M$ is defined by
  \[  \sv M := \lone{\fcl M _\R}
            = \inf \bigl\{ \lone c
                   \bigm|  \text{$c \in \csing * M$ is an
                                 $\R$-fundamental cycle of~$M$}
                   \bigr\},
  \]
  where $\fcl M _\R \in \hsing * M$ denotes the image of the
  fundamental class~$\fcl M \in \hsing * {M;\Z}$ of~$M$ under the
  change of coefficients homomorphism~$\hsing * {M;\Z} \longrightarrow
  \hsing * M$.
  
  The simplicial volume of spheres and tori is zero; on the other
  hand, the simplicial volume of oriented, closed, hyperbolic
  surfaces~$S$ equals~$2\cdot|\chi(S)|$~\cite[p.~8f]{vbc}, and more
  generally the simplicial volume of closed Riemannian manifolds of
  negative curvature is non-zero~\cite{inoue,thurston}.

  The simplicial volume of oriented, closed, connected, locally
  symmetric spaces of non-compact type is
  non-zero~\cite{lafontschmidt}. Therefore, if it were known that the
  fundamental groups of oriented, closed, connected, irreducible
  locally symmetric spaces of non-compact type have amenable
  centralisers, one could apply Theorem~\ref{svthm} also in this
  case.

  As an indication of the applicability of Theorems~\ref{asphthm} and
  Theorem~\ref{svthm}, we mention the class of word-hyperbolic groups
  (cf.~Section~\ref{hypgroupsubsec}). 

\subsection*{Organisation} 

This article is organised as follows: In Section~\ref{hypproofsec}, we
derive Theorem~\ref{hypthm} from the more general statements in
Theorem~\ref{asphthm} and Theorem~\ref{svthm}. Theorems~\ref{asphthm}
and~\ref{svthm} in turn are proved in Section~\ref{proofsec}.

\subsection*{Acknowledgements}

I would like to thank Wolfgang L\"uck for clarifying discussions.

%
\section{Deriving the negatively curved case}\label{hypproofsec}

\noindent
To show that fundamental classes of negatively curved manifolds cannot
be represented by products of manifolds (Theorem~\ref{hypthm}), it
suffices to establish this statement as a special case of
Theorem~\ref{svthm}:

\subsection{Word-hyperbolic groups}\label{hypgroupsubsec}

The notion of word-hyperbolicity is a group-theoretic analogue of
negative curvature introduced by Gromov~\cite{gromovhyperbolic}. For
example, fundamental groups of oriented, closed, connected Riemannian
manifolds of negative curvature are word-hyperbolic by the \v
Svarc-Milnor lemma~\cite[Section~V.D]{delaharpe}.

\begin{bem}\label{hyprem}\hfil
  \begin{itemize}
    \item Word-hyperbolic groups have amenable centralisers: Let $G$
      be a word-hyperbolic group. If $g \in G$ has
      infinite order, then the subgroup of~$G$ generated by~$g$ has
      finite index in the centraliser~$C_G(g)$~\cite[Corollary~3.10
      on~p.~462]{bh}. Therefore, the group~$C_G(g)$ contains a normal, 
      infinite cyclic subgroup of finite index and hence is
      amenable~\cite{paterson}.  
    \item All infinite word-hyperbolic groups are large in the sense of
      Definition~\ref{largedef} above~\cite[Proposition~2.22
      on~p.~458]{bh}. 
      \qedhere
  \end{itemize}
\end{bem}

\subsection{Proof of Theorem~\ref{hypthm}}

Using the facts on word-hyperbolic groups listed above, we can show
that fundamental classes of negatively curved Riemannian manifolds
cannot be represented by products of manifolds:

\begin{bew}[Proof (of Theorem~\ref{hypthm})]
   Because $M$ is an oriented, closed, connected, smooth manifold
   admitting a Riemannian metric of negative sectional curvature,
   Thurston's straightening technique shows that $M$ has non-zero
   simplicial volume~\cite{inoue, thurston}.

   By the Cartan-Hadamard theorem, $M$ is aspherical; in particular,
   $\pi_1(M)$ is torsion-free and we may assume that $\pi_1(M)$ is
   non-trivial; hence, $\pi_1(M)$ is large.
   Moreover, the fundamental group of~$M$ is word-hyperbolic by the \v
   Svarc-Milnor Lemma.  Therefore, Remark~\ref{hyprem} shows that
   $\pi_1(M)$ has amenable centralisers.

   Applying Theorem~\ref{svthm} to~$M$ yields that the fundamental
   class of~$M$ cannot be represented by a product of manifolds.
\end{bew}

%
\section{Proof of Theorems~\ref{asphthm} and~\ref{svthm}}\label{proofsec}

\noindent
In this section, we prove Theorems~\ref{asphthm}
and~\ref{svthm}. After introducing the relevant notation, we collect
in Section~\ref{groupscenesubsec} a number of observations on the
obstructions to representability by products on the level of
fundamental groups. In Section~\ref{concludesubsec}, we conclude the
proofs of Theorem~\ref{asphthm} and~\ref{svthm}.

\subsection{Notation}\label{notationsubsec}

For the remainder of this section, we assume that $M$ is an oriented,
closed, connected manifold whose fundamental class of~$M$ can be
represented by a product of manifolds and whose fundamental group has
centralisers in a class~$C$ with products. 
Combining multiple factors we may assume that the fundamental class
of~$M$ can be represented by a product of two manifolds, i.e., there
are~$d \in \Z\setminus\{0\}$, and oriented, closed, connected
manifolds~$S_1$,~$S_2$ of non-zero dimension admitting a continuous
map~$f \colon S_1 \times S_2 \longrightarrow M$ such that
\[ \hsing n {f;\Z}(\fcl{S_1 \times S_2}) = d \cdot \fcl M,
\]
where $n := \dim M$; in particular, $\dim S_1 + \dim S_2 = n$ and
$\deg f = d$. Notice that $n = \dim M > 0$ because the fundamental
class of~$M$ can be represented by a product of manifolds.

In the course of the proofs, we use the following
abbreviations:

    \begin{itemize}
      \item We choose base-points~$s_1 \in S_1$ and $s_2 \in S_2$, and
        set~$m := f(s_1, s_2) \in M$. In the following, all
        fundamental groups are taken with respect to these base-points
        or base-points deduced from these in an obvious way.
      \item We write
        \begin{align*}
          f_1 & := f|_{S_1 \times \{s_2\}} 
                   \colon S_1 \times \{s_2\} \longrightarrow M,\\
          f_2 & := f|_{\{s_1\} \times S_2} 
                   \colon \{s_1\} \times S_2 \longrightarrow M.
        \end{align*}
      \item 
        Correspondingly, on the level of fundamental groups, we define
        \begin{align*}
          H_1 & := \im \pi_1(f_1) \subset \pi_1(M,m),\\
          H_2 & := \im \pi_1(f_2) \subset \pi_1(M,m),
        \end{align*}
        and
        \begin{align*}
          H_1' & := \pi_1(S_1 \times \{s_2\}, (s_1,s_2)) 
                    \subset \pi_1(S_1 \times S_2, (s_1,s_2)),\\
          H_2' & := \pi_1(\{s_1\} \times S_2, (s_1,s_2)) 
                    \subset \pi_1(S_1 \times S_2, (s_1,s_2)).
        \end{align*}
      \item Finally, $G := \im \pi_1(f) \subset \pi_1(M,m)$.
        \qedhere
    \end{itemize}

\subsection{The scene on the level of fundamental groups}\label{groupscenesubsec}

As indicated in the introduction, we now analyse the corresponding
scene on the level of fundamental groups: 

  \begin{enumerate}
    \item\label{centitem}  
      \emph{The group~$G$ has centralisers in~$C$ and $G$ has finite
        index in~$\pi_1(M)$.}
      \begin{bew}
      Centralisers in~$G$ are subgroups of centralisers
      in~$\pi_1(M)$, and hence are subgroups of groups in~$C$. Because
      $C$ is closed with respect to taking subgroups, the group~$G$ has
      centralisers in~$C$.

      Furthermore, $G$ has finite index in~$\pi_1(M)$: Let $p \colon
      \overline M \longrightarrow M$ be the covering associated with
      the subgroup~$G \subset \pi_1(M)$; in particular, $\overline M$
      is a manifold with~$\dim \overline M = \dim M$. By definition
      of~$G$, covering theory provides us with a lift~$\overline f
      \colon S_1 \times S_2 \longrightarrow \overline M$ of~$f$, i.e.,
      $p \circ \overline f = f$.
      Because $\deg f \neq 0$, applying top homology shows that
      $\overline M$ must be compact and $\deg f = \deg p \cdot \deg
      \overline f$. Therefore, the index
      \[ [ \pi_1(M) : G ] 
         = \left|\deg p\right|
      \]
      must be finite.
      \end{bew}
    \item\label{commgenitem} \emph{The set~$H_1 \cup H_2$ generates~$G$,
        and $H_1$ and $H_2$ commute with each other.}
      \begin{bew}    
      The inclusions~$H_1' \hookrightarrow \pi_1(S_1
      \times S_2)$ and~$H_2' \hookrightarrow \pi_1(S_1
      \times S_2)$ induce an isomorphism~$H_1' \times H_2' \cong
      \pi_1(S_1 \times S_2)$. In particular, $\pi_1(S_1 \times S_2)$
      is generated by~$H_1' \cup H_2'$. Therefore
      \[ H_1 \cup H_2 = \pi_1(f) (H_1' \cup H_2')
      \]
      generates~$\im \pi_1(f)$, which coincides -- by definition --
      with~$G$. 

      Because~$H_1'$ and~$H_2'$ commute with each other in~$\pi_1(S_1
      \times S_2) \cong H_1' \times H_2'$, it follows that also their
      images~$H_1$ and~$H_2$ under the homomorphism~$\pi_1(f)$
      commute.
      \end{bew}
    \item\label{quotitem}
      \emph{The group~$G$ is a quotient of~$H_1 \times H_2$.} 
      \begin{bew}
        Because $H_1$ and~$H_2$ commute, there is a homomorphism~$H_1
        \times H_2 \longrightarrow G$, which is surjective by
        part~\ref{commgenitem}. Hence, $G$ is a quotient of~$H_1 \times
        H_2$.
      \end{bew}
    \item\label{largeitem} 
      \emph{If $H_1$ is large (in the sense of
        Definition~\ref{largedef}), then $H_2$ lies
        in~$C$, and vice versa .}
      \begin{bew}
        If $H_1$ is large, then there is an element~$g \in H_1$ of
        infinite order, and hence part~\ref{centitem} implies
        that~$C_G(g) \in C$. On the other hand, $H_2 \subset C_G(g)$ by
        part~\ref{commgenitem}.
      \end{bew}
    \item\label{largeritem}
      \emph{If $\pi_1(M)$ is large, then at least one of the
        groups~$H_1$ and~$H_2$ is large.} 
      \begin{bew}
        In view of part~\ref{quotitem} it suffices to show that $G$ is
        large:
        If~$\pi_1(M)$ is large, we find an element~$g \in
        \pi_1(M)$ of infinite order. Because $G$ has finite index
        in~$\pi_1(M)$, there are~$m$,~$n \in \N_{>1}$ with $g^m \cdot G
        = g^n \cdot G$ and $m \neq n$. Hence, $g^{m-n} \in G$. On the
        other hand, $m-n \neq 0$ implies that $g^{m-n}$ has infinite
        order. I.e., $G$ is large.
      \end{bew}
    \item\label{diagitem}
      \emph{The diagram in Figure~\ref{classfig} is commutative. In particular,
        there are homology classes~$\alpha_1 \in \hsing{\dim 
          S_1}{BH_1}$ and~$\alpha_2 \in \hsing{\dim S_2}{BH_2}$ 
        satisfying
        \[ H_n(c_M)(d \cdot \fcl M_\R) = H_n(\overline\varphi)(\alpha_1 \times
                                                       \alpha_2). 
        \]}
      \begin{bew}
      \begin{figure}
        \begin{align*}
          \xymatrix@C=3em{%
                S_1 \times S_2 
                \ar[rrrr]^-f
                \ar[dr]^{c_{S_1} \times c_{S_2}}
                \ar[ddd]_-{c_{S_1 \times S_2}}
           &&&& M 
                \ar[ddd]^-{c_M} 
                \\
              & BH_1' \times BH_2'
                \ar[rr]^-{B\pi_1(f_1) \times B\pi_1(f_2)}
             && BH_1 \times BH_2 
                \ar@{-->}[ddr]^{\overline \varphi}
                \\
              & B(H_1' \times H_2')
                \ar@{<->}[u]^-{\simeq}
                \ar[rr]_-{B(\pi_1(f_1) \times \pi_1(f_2))}
                \ar[dl]^-{\quad\smash{\varphi'}}
             && B(H_1 \times H_2)
                \ar@{<->}[u]_-{\simeq}
                \ar[dr]_-{\smash{\varphi}\quad}
                \\
                B\pi_1(S_1 \times S_2) 
                \ar[rrrr]_-{B\pi_1(f)}
           &&&& B\pi_1(M)
            }
        \end{align*}
      \caption{Proof of part~\ref{diagitem}}
      \label{classfig}
      \end{figure}
        We first explain the notation occurring in
        Figure~\ref{classfig}: For a path-connected space~$X$, we
        write~$c_X \colon X \longrightarrow B\pi_1(X)$ for the
        classifying map. Recall that every homomorphism~$\psi \colon
        K' \longrightarrow K$ of groups yields a continuous map~$B\psi
        \colon BK' \longrightarrow BK$ that induces the given
        homomorphism~$\psi$ on the level of fundamental groups;
        moreover, $B\psi$ is characterised uniquely up to homotopy by
        this property.

        The vertical homotopy equivalences in the centre of the
        diagram are induced by the projections/inclusions on the level
        of groups.  The map~$\varphi'$ is induced by the canonical
        isomorphism~$H_1' \times H_2' \longrightarrow \pi_1(S_1 \times
        S_2)$ given by the inclusions. Finally, for~$\varphi$ we
        observe that the inclusions~$H_1 \longrightarrow \pi_1(M)$
        and~$H_2 \longrightarrow \pi_1(M)$ give rise to a
        homomorphism~$H_1 \times H_2 \longrightarrow \pi_1(M)$,
        because $H_1$ and $H_2$ commute in~$G \subset \pi_1(M)$; the
        map~$\varphi$ is the continuous map induced by this
        homomorphism.
      
        It is a routine matter to verify that the diagram in
        Figure~\ref{classfig} is commutative up to homotopy. Hence,
        the corresponding diagram on the level of singular homology
        with $\R$-coefficients is commutative.

        In the following, we abbreviate the composition of~$\varphi$
        with the homotopy equivalence~$BH_1 \times BH_2
        \longrightarrow B(H_1 \times H_2)$ by $\overline \varphi
        \colon BH_1 \times BH_2 \longrightarrow B\pi_1(M).$ Using
        the naturality of the homological cross-product, we obtain in
        singular homology with $\R$-co\-ef\-fi\-cients the relation
        \begin{align*}
            \hsing n {c_M} (d \cdot \fcl M _\R)
        & = \hsing n {c_M} \circ \hsing n {f}(\fcl{S_1 \times S_2}_\R)\\
        & = \hsing n {c_M} \circ \hsing n {f}(\fcl{S_1}_\R \times \fcl{S_2}_\R)\\
        & =       \hsing n {\overline \varphi}
                  \circ \hsing[big] n {B\pi_1(f_1) \times B\pi_1(f_2) \circ c_{S_1} \times c_{S_2}}
                        (\fcl{S_1\times S_2}_\R) \\
        & =       \hsing n {\overline \varphi}(\alpha_1 \times \alpha_2),
        \end{align*}
        where we put $\alpha_j :=  \hsing {\dim S_j} {B\pi_1(f_j) \circ
          c_{S_j}}(\fcl{S_j}_\R)$ for~$j \in \{1,2\}$.
      \end{bew}
    \end{enumerate}

\subsection{Completing the proofs}\label{concludesubsec}

Using the facts on the groups~$H_1$ and~$H_2$ established in the
previous section, we can finish the proofs of Theorem~\ref{asphthm}
and~\ref{svthm}:

\begin{bew}[Proof (of Theorem~\ref{asphthm})]
  By hypothesis, the manifold~$M$ is aspherical; in particular,
  $\pi_1(M)$ is torsion-free and non-trivial (because $\dim M >0$) and
  hence any non-trivial subgroup of~$\pi_1(M)$ is large.  

  If $H_1$ were trivial, then $\alpha_1 = 0$ in part~\ref{diagitem},
  which would show that $H_n(c_M)(\fcl M_\R) = 0$. However, $M$ is
  aspherical and therefore, $c_M$ is a homotopy
  equivalence. Analogously, $H_2$ cannot be trivial; i.e., $H_1$ and
  $H_2$ are large.  

  Therefore, we obtain from part~\ref{largeitem} that $H_1$ and $H_2$
  lie in~$C$. Using part~\ref{quotitem} and the fact that $C$ is
  closed with respect to taking finite products, subgroups, and
  quotients, we deduce that $G$ is in~$C$. Furthermore, $G$ has finite
  index in~$\pi_1(M)$ by part~\ref{centitem}.
\end{bew}

\begin{bew}[Proof (of Theorem~\ref{svthm})]
  By hypothesis, $\pi_1(M)$ is large and $\sv M > 0$. \emph{Assume}
  that the fundamental class of~$M$ can be represented by a product of
  manifolds; thus, the observations of Section~\ref{groupscenesubsec}
  apply.

  By part~\ref{largeritem}, we may assume that $H_1$ is large and
  hence that $H_2$ is amenable (part~\ref{largeitem}).

  We now apply the $\ell^1$-semi-norm to the equation of
  part~\ref{diagitem}. At this point, we rely on the following classic
  results on the $\ell^1$-semi-norm:
      \begin{itemize}
        \item The homomorphism~$\hsing n {c_M}
          \colon \hsing n M \longrightarrow \hsing n {B\pi_1(M)}$ is
          isometric with respect to the
          $\ell^1$-semi-norm by the mapping theorem in bounded
          cohomology~\cite[p.~40/18, Theorem~4.3]{vbc, ivanov}.  
        \item The $\ell^1$-semi-norm on connected, countable
          \cw-complexes with amenable fundamental group is zero in
          non-zero degree~\cite[p.~40/18, Theorem~4.3]{vbc, ivanov}.
        \item The standard triangulation of products of simplices
          shows that the $\ell^1$-semi-norm is compatible with the
          homological cross-product in the following
          sense: For all spaces~$X$,~$Y$ and all~$\alpha \in \hsing k
          X$,~$\beta \in \hsing \ell Y$, we 
          have
          \[ \lone{\alpha \times \beta} 
             \leq 2^{k+\ell} \cdot \lone \alpha \cdot \lone \beta.
          \]
        \item
          The $\ell^1$-semi-norm is a functorial semi-norm, i.e.,
          the map~$\hsing n {\overline \varphi}$ does not increase
          the $\ell^1$-semi-norm.
      \end{itemize}
      
      Since $H_2$ is countable, there exists a countable \cw-model
      of~$BH_2$.  Because $H_2$ is amenable, the properties
      of~$\lone{\args}$ listed above and part~\ref{diagitem} imply
      that
      \begin{align*} 
                 \sv M
        & =      \frac{1}{|d|} \cdot \lone[big]{\hsing n {c_M}
                                            (d \cdot \fcl M_\R)} 
          =      \frac 1{|d|} \cdot  \lone[big]{\hsing n {\overline
                                        \varphi}(\alpha_1 \times
                                      \alpha_2)} \\
        & \leq         \frac{2^n}{|d|} \cdot \lone{\alpha_1} \cdot
        \lone{\alpha_2} \\
        & = 0,
      \end{align*}
      which contradicts the hypothesis that the simplicial volume~$\sv
      M$ of~$M$ is non-zero. 
\end{bew}


\vfill

\noindent
\begin{minipage}{\linewidth}
\small
\noindent
\emph{Clara L\"oh} \\[1ex]
\begin{tabular}{@{}ll}
  address:       & Fachbereich Mathematik und Informatik der
                   \wwu~M\"unster\\
                 & Einsteinstr.~62\\
                 & 48149~M\"unster, Germany \\[1ex]
  email:         & \textsf{clara.loeh@uni-muenster.de}\\
  \smaller{URL}: & \textsf{http://wwwmath.uni-muenster.de/u/clara.loeh}
\end{tabular}
\end{minipage}

\begin{thebibliography}{99}
  \small
  \setlength{\itemsep}{0pt}
 \bibitem{bh}M.R.~Bridson, A.~Haefliger. 
          \emph{Metric Spaces of Non-positive Curvature},
          Volume~319 of \emph{Grundlehren der Mathematischen
          Wissenschaften}, 
          Springer, 1999.

 \bibitem{vbc} M.~Gromov. Volume and bounded cohomology. 
          \emph{Publ.\ Math.\ IHES}, 56, pp.~5--99, 1982.%

 \bibitem{gromovhyperbolic} M.~Gromov. \emph{Hyperbolic groups}. 
          \emph{Essays in Group Theory} (edited by S.M.~Gersten),
           Math. Sci. Res. Inst. Publ., 8, pp.~75--263, Springer, 1987. 

 \bibitem{gromovmetric} M.~Gromov. \emph{Metric Structures for
          Riemannian and Non-Riemannian Spaces} with appendices by
          M.~Katz, P.~Pansu and S.~Semmes, translated from the French
          by Sean Michael Bates. Volume~152 of  
          \emph{Progress in Mathematics}, Birkh\"auser, 1999.

 \bibitem{delaharpe} P.~de~la~Harpe.
          \emph{Topics in Geometric Group Theory}, 
          Chicago University Press, 2000.

  \bibitem{inoue} H.~Inoue, K. Yano. The Gromov invariant
          of negatively curved manifolds. \emph{Topology}, 21,
          No.~1, pp.~83--89, 1981.

 \bibitem{ivanov} N.V.~Ivanov. Foundations of the theory of
          bounded cohomology. \emph{J.~Soviet Math.}, 37,
          pp.~1090--1114, 1987.

   \bibitem{lafontschmidt} J.-F.~Lafont, B.~Schmidt. Simplicial volume of
          closed locally symmetric spaces of non-compact
          type. \emph{Acta Math.}, 197, No.~1, pp.~129--143, 2006.  

 \bibitem{paterson} A.L.T.~Paterson. \emph{Amenability}. Volume~29 of
          \emph{Mathematical Surveys and Monographs}, AMS, 1988.

 \bibitem{thom} M.R.~Thom. Sur un probl\`eme de
          Steenrod. \emph{C.~R.~Acad.~Sci.~Paris}, 236,
          pp.~1128--1130, 1953.

 \bibitem{thurston} W.P.~Thurston. \emph{Geometry and
          Topology of 3-Manifolds}. Lecture notes, Princeton, 1978. 
          Available at \url{http://www.msri.org/publications/books/gt3m}.

\end{thebibliography}
\end{document}